\documentclass[10pt]{article}
\usepackage{amsmath}
\usepackage{amsfonts}
\usepackage{amssymb}
\long\def\proof#1{\removelastskip\vskip\baselineskip\relax\noindent{\it
Proof\if!#1!\else\ \ignorespaces#1\fi.\ }\ignorespaces}

\newcommand{\lgs}[2]{\mbox{$\left(\frac{#1}{#2}\right)$}}
\newcommand{\leg}[2]{\mbox{$\left(\dfrac{#1}{#2}\right)$}}

\newcommand{\ov}[1]{\overline{\vphantom{T}#1}}
\newcommand{\Os}{\tilde{O}}
\newcommand{\Q}{{\mathbb Q}}
\newcommand{\Z}{{\mathbb Z}}

\newcommand{\C}{{\mathbb C}}

\newcommand{\z}{\zeta}
\renewcommand{\th}{\theta}
\renewcommand{\H}{{\mathcal H}}
\newcommand{\la}{\lambda}

\newcommand{\D}{\Delta}
\newcommand{\G}{\Gamma}

\newcommand{\ga}{\gamma}

\newcommand{\dd}{\diagdown}

\newcommand{\eps}{\varepsilon}

\newcommand{\psmm}[4]{\left(\begin{smallmatrix}{#1}&{#2}\\{#3}&{#4}\end{smallmatrix}\right)}
\newcommand{\g}{{\mathfrak g}}

\newcommand{\Proof}{{\it Proof. \/}}
\newcommand{\squareforqed}{\hbox{\rlap{$\sqcap$}$\sqcup$}}
\newcommand{\qed}{\ifmmode\squareforqed\else{\unskip\nobreak\hfil
\penalty50\hskip1em\null\nobreak\hfil\squareforqed
\parfillskip=0pt\finalhyphendemerits=0\endgraf}\fi}

\newcommand{\fp}{\qed\removelastskip\vskip\baselineskip\relax}

\newtheorem{theorem}{Theorem}[section]
\newtheorem{corollary}[theorem]{Corollary}
\newtheorem{proposition}[theorem]{Proposition}

\newtheorem{definition}[theorem]{Definition}
\newtheorem{conjecture}[theorem]{Conjecture}

\newcommand{\litem}{\par\noindent\dimen0=\parindent%
    \advance\dimen0 by-4pt
               \hangindent=\dimen0\ltextindent}

\newcommand{\ltextindent}[1]{\hbox to \hangindent{#1\hss}\ignorespaces}
\newcommand{\ltextjndent}[1]{\hbox to \hangindent{#1\hss}\ignorespaces\kern-1ex}

\renewcommand{\pmod}[1]{\allowbreak\ ({\rm{mod}}\,\,#1)}

\begin{document}
\pagestyle{plain}

\title{Computing $L$-Functions of Quadratic Characters at Negative Integers}
\author{Henri Cohen\\
Universit\'e de Bordeaux,\\
Institut de Math\'ematiques, U.M.R. 5251 du C.N.R.S,\\
351 Cours de la Lib\'eration,\\
33405 TALENCE Cedex, FRANCE}

\maketitle
\begin{abstract}
  We survey a number of different methods for computing $L(\chi,1-k)$ for
  a Dirichlet character $\chi$, with particular emphasis on quadratic
  characters. The main conclusion is that when $k$ is not too large
  (for instance $k\le100$) the best method comes from the use of Eisenstein
  series of half-integral weight, while when $k$ is large the best method
  is the use of the complete functional equation, unless the conductor of
  $\chi$ is really large, in which case the previous method again prevails.
\end{abstract}

\smallskip

\section{Introduction}

This paper can be considered as a complement of two of my old papers
\cite{CohL} and \cite{CohH}, updated to include new formulas, and surveying
existing methods.

The general goal of this paper is to give efficient methods for computing
the values at negative integers of $L$-functions of Dirichlet characters
$\chi$. Since these values are algebraic numbers, more precisely belong to the
cyclotomic field $\Q(\chi)$, we want to know their \emph{exact} value.
When $\chi(-1)=(-1)^{r-1}$ we have $L(\chi,1-k)=0$, so we always assume
implicitly that $\chi(-1)=(-1)^r$. In addition, if $\chi$ is a non-primitive
character modulo $F$ and $\chi_f$ is the primitive character associated to
$\chi$, we have
$$L(\chi,1-k)=L(\chi_f,1-k)\prod_{p\mid F,\ p\nmid f}(1-\chi_f(p)p^{r-1})\;,$$
so we may assume that $\chi$ is primitive.

Note that we will not consider the slightly different problem of computing
\emph{tables} of $L(\chi,1-k)$, either for fixed $k$ and varying $\chi$
(such as $\chi=\chi_D$ the quadratic character of discriminant $D$), or
for fixed $\chi$ and varying $k$, although several of the methods considered
here can be used for this purpose.

\medskip

In addition to their intrinsic interest, these computations have several
applications, for instance:

\begin{enumerate}\item Computing $\la$-invariants of quadratic fields
  (I am indebted to J.~Ellenberg and S.~Jain for this, see \cite{EJV}).
\item Computing Sato--Tate distributions for modular forms of
  half-integral weight, see \cite{IW} and \cite{IOTW}.
\item Computing Hardy--Littlewood constants of polynomials, see \cite{Bel-Coh}.
\end{enumerate}

\medskip

There exist at least five different methods for computing these quantities,
some having several variants. We denote by $F$ the conductor of $\chi$.

\begin{enumerate}
\item Bernoulli methods: one can express $L(\chi,1-k)$ as a finite sum
involving $O(F)$ terms and Bernoulli numbers, so that the time required is
$\Os(F)$ (we use the ``soft-O'' notation $\Os(X)$ to mean $O(X^{1+\eps})$ for
any $\eps>0$). This method has two variants: one which uses directly the
definition of $\chi$-Bernoulli numbers, the second which uses
\emph{recursions}.
\item Use of the \emph{complete functional equation}. Using it, it is
sufficient first to compute numerically $L(\ov{\chi},k)$ to sufficient accuracy
(given by the functional equation), which is done using the Euler product,
and second to know an upper bound on the denominator of $L(\chi,1-k)$,
which is easy (and usually equal to $1$). The required time is also $\Os(F)$,
but with a much smaller implicit $O()$ constant.
\item Use of the \emph{approximate functional equation}, which involves in
particular computing the incomplete gamma function or similar higher
transcendental functions. The required time is $\Os(F^{1/2})$, but with a
large implicit $O()$ constant.
\item Use of Hecke-Eisenstein series (Hilbert modular forms) on the full
modular group, which expresses $L(\chi,1-k)$ as a finite sum involving
$O(F^{1/2})$ terms and (twisted) sum of divisors function. The required time
is $\Os(F^{1/2})$ with a very small implicit $O()$ constant. A variant
which is useful only for very small $k$ such as $k\le10$ uses Hecke-Eisenstein
series on congruence subgroups of small level.
\item Use of Eisenstein series of half-integral weight over $\G_0(4)$, which
again expresses $L(\chi,1-k)$ as a finite sum involving $O(F^{1/2})$ terms
and (twisted) sum of divisors function, but different from the previous ones.
The required time is again $\Os(F^{1/2})$, but with an even smaller implicit
$O()$ constant. An important variant, valid for all $k$, is to use
modular forms of half-integral weight on subgroups of $\G_0(4)$.
\end{enumerate}

The first three methods are completely general, but the last two are really
efficient only if $\chi$ is equal to a quadratic character or possibly a
quadratic character times a character of small conductor. We will
therefore present all five methods and their variants, but consider the
last two methods only in the case of quadratic characters, and therefore
compare them only in this case.

\medskip

After implementing these methods and comparing their running times for various
values of $F$, we have arrived at the following conclusions: first,
the two fastest methods are always either the fifth (Eisenstein series
of half-integral weight) or the second (complete functional equation).
Second, one should choose the second method only if $k$ is large, for instance
$k\ge100$, except if $F$ is large. Note also that the case $F=1$ corresponds to
the computation of Bernoulli numbers, and that indeed the fastest method
for this is the use of the complete functional equation of the Riemann
zeta function.

Because of these conclusions, we will give explicitly the formulas for
the first, third, and fourth method, but only describe the precise
implementations and timings for the second and fifth, which are the really
useful ones.

\section{Bernoulli Methods}

\subsection{Direct Formulas}

\begin{proposition}
Define the $\chi$-Bernoulli numbers $B_k(\chi)$ by the  generating function
$$\dfrac{T}{e^{FT}-1}\sum_{0\le r<F}\chi(r)e^{rT}
=\sum_{k\ge0}\dfrac{B_k(\chi)}{k!}T^k\;.$$
Then
$$L(\chi,1-k)=-\dfrac{B_k(\chi)}{k}-\chi(0)\delta_{k,1}\;.$$
\end{proposition}

Note that since we assume $\chi$ primitive, the term 
$\chi(0)\delta_{k,1}$ vanishes unless $F=1$ and $k=1$, in which case
$L(\chi,1-k)=\z(0)=-1/2$. Also, recall that for $k\ge2$ we have
$B_k(\chi)=0$ if $\chi(-1)\ne(-1)^k$.

\begin{proposition}\label{prop1} Set $S_n(\chi)=\sum_{0\le r<F}\chi(r)r^n$.
We have
\begin{align*}
B_k(\chi)&=\dfrac{1}{F}\left(S_k(\chi)-\dfrac{kF}{2}S_{k-1}(\chi)
+\sum_{1\le j\le k/2}\binom{k}{2j}B_{2j}F^{2j}S_{k-2j}(\chi)\right)\\
&=\dfrac{1}{F}\sum_{0\le r<F}\chi(r)\left(r^k-\dfrac{kF}{2}r^{k-1}+\sum_{1\le j\le k/2}\binom{k}{2j}B_{2j}r^{k-2j}F^{2j}\right)\\
&=\dfrac{1}{F}\sum_{1\le j\le k+1}\dfrac{(-1)^{j-1}}{j}\binom{k+1}{j}\sum_{0\le r<Fj}\chi(r)r^k\;.\end{align*}
\end{proposition}

\subsection{Recursions}

There are a large number of recursions for $B_k(\chi)$. The following
three propositions give some of the most important ones:

\begin{proposition}
We have the recursion
$$\sum_{0\le j<k}F^{k-j}\binom{k}{j}B_j(\chi)=kS_{k-1}(\chi)\;,$$
where $S_n(\chi)$ is as above.
\end{proposition}

\begin{proposition}\label{casse2} Let $\chi$ be a nontrivial primitive 
character of conductor $F$, set $\eps=\ov{\chi}(2)$ and
$$Q_k(\chi)=\sum_{1\le r<F/2}\chi(r)r^k\;.$$
We have the recursion
$$(2^k-\eps)B_k(\chi)=-\Biggl(k2^{k-1}Q_{k-1}(\chi)+\sum_{1\le j<k/2}\binom{k}{2j}(2^{k-1-2j}-\eps)F^{2j}B_{k-2j}(\chi)\Biggr)\;.$$
\end{proposition}

\begin{proposition}\label{lucas} Let $\chi$ be a nontrivial primitive 
character of conductor $F$.
\begin{enumerate}\item If $\chi$ is even we have
  $$\sum_{0\le j\le (k-1)/2}\binom{k}{2j+1}F^{2j}\dfrac{B_{2k-2j}(\chi)}{2k-2j}=\dfrac{(-1)^k}{F}\sum_{0\le r<F/2}\chi(r)r^k(F-r)^k\;.$$
\item If $\chi$ is odd we have
\begin{align*}\sum_{0\le j\le (k-1)/2}\binom{k}{2j+1}&F^{2j}B_{2k-1-2j}(\chi)=\\
&=\dfrac{(-1)^kk}{F}\sum_{0\le r<F/2}\chi(r)r^{k-1}(F-r)^{k-1}(F-2r)\;.\end{align*}
\end{enumerate}\end{proposition}

In practice, it seems that the fastest way to compute $L(\chi,1-k)$ using
$\chi$-Bernoulli numbers is to use Proposition \ref{casse2}, but
it is not competitive with the other methods that we are going to give.

\section{Using the Complete Functional Equation}

In this section and the next, we use approximate methods to compute
$L(\chi,1-k)$, which for simplicity we call \emph{transcendental} methods,
since they use transcendental functions.

Since our goal is to compute these values as exact algebraic numbers,
and since we know that $L(\chi,1-k)\in\Q(\z_u)$, where $u$ is the order of
$\chi$, we simply need to know an upper bound for the denominator of
$L(\chi,1-k)$ as an algebraic number, and we need to compute simultaneously
$L(\chi^j,1-k)$ for $k$ modulo $u$ and coprime to $u$, so that the individual
values can then be obtained by simple linear algebra. A priori this involves
$\phi(u)$ computations, but since $L(\chi^{-1},1-k)$ is simply the complex
conjugate of $L(\chi,1-k)$, only $\lceil\phi(u)/2\rceil$ computations are
needed. In particular, if $u=1$, $2$, $3$, $4$, or $6$, a single computation
suffices.

\smallskip

Thus, we need two types of results: one giving the approximate size of
$L(\chi,1-k)$, so as to determine the relative accuracy with which to do
the computations, and second an upper bound for its denominator. The first
result is standard, and the second can be found in Section 11.4 of \cite{Coh2}:

\begin{proposition} We have
$$L(\chi,1-k)=\dfrac{2\cdot(k-1)!F^k}{(-2i\pi)^k\g(\ov{\chi})}L(\ov{\chi},k)\;,$$
where $\g(\ov{\chi})$ is the standard Gauss sum of modulus $|F|^{1/2}$
associated to $\ov{\chi}$.\end{proposition}

\begin{corollary} As $k\to\infty$ we have
$$|L(\chi,1-k)|\sim 2\cdot e^{-1/2}\left(\dfrac{kF}{2\pi e}\right)^{k-1/2}\;.$$
\end{corollary}

\Proof Clear from Stirling's formula and the fact that $L(\ov{\chi},k)$ tends
to $1$ when $k\to\infty$.\fp

\begin{theorem} Denote by $u$ the order of $\chi$, so that $u\mid\phi(F)$
and $L(\chi,1-k)\in K=\Q(\z_u)$. We have $D(\chi,k)L(\chi,1-k)\in\Z[\z_u]$,
where the ``denominator'' $D(\chi,k)$ can be chosen as follows:
\begin{enumerate}\item If $F$ is not a prime power then $D(\chi,k)=1$.
\item Assume that $F=p^v$ for some odd prime $p$ and $v\ge1$.
\begin{enumerate}\item If $u\ne p^{v-1}(p-1)/\gcd(p-1,k)$ then $D(\chi,k)=1$.
\item If $u=p^{v-1}(p-1)/\gcd(p-1,k)$ then $D(\chi,k)=pk/((p-1)/u)$ if $v=1$ or
$D(\chi,k)=\chi(1+p)-1$ if $v\ge2$.
\end{enumerate}
\item If $F=2^v$ for some $v\ge2$ then $D(\chi,k)=1$ if $v\ge3$, while
$D(\chi,k)=2$ if $v=2$.
\item If $F=1$ then $D(\chi,k)=k\prod_{(p-1)\mid k}p$.
\end{enumerate}\end{theorem}

Stronger statements are easy to obtain, see \cite{Coh2}, but these bounds
are sufficient.

\smallskip

To compute $L(\chi,1-k)$ using these results, we proceed as follows.
Let $B$ be chosen so that 
$$B>(k-1/2)\log(kF/(2\pi e))+\log(|D(\chi,k)|)+10\;,$$
where $10$ is simply a safety margin. Thanks to the above two results,
computing $L(\chi,1-k)$ to \emph{relative} accuracy $e^{-B}$ will guarantee
that the coefficients of the algebraic integer $D(\chi,k)L(\chi,1-k)$ on
the integral basis $(\z_u^j)_{0\le j<\phi(u)}$ will be correct to accuracy
$e^{-5}$, say, and since they are in $\Z$, they can thus be recovered
exactly.

Thanks to the functional equation, it is thus sufficient to compute
$L(\ov{\chi},k)$ to relative accuracy $e^{-B}$, but since $L(\ov{\chi},k)$
is close to $1$, $k$ being large, this is the same as absolute accuracy.
Note from the above formula that $B$ will be (considerably) larger than $k$.

To compute $L(\ov{\chi},k)$, we first compute 
$\prod_{p\le L(B,k)}(1-\chi(p)/p^k)$, using an internal accuracy of 
$e^{-kB/(k-1)}$ and limit $L(B,k)=(e^B/(k-1))^{1/(k-1)}$.
More precisely, we initially set $P=1$, and for primes $p$ going from
$2$ to $L(B,k)$, we compute $1/p^k$ to $p^ke^{-kB/(k-1)}$ of 
relative accuracy (this is crucial), and then set $P\gets P-P(1/p^k)$.
It is clear that this will compute $1/L(\chi,k)$ to the desired precision,
from which we immediately obtain $L(\ov{\chi},k)$. Important implementation
remark: to compute the accuracy needed in the intermediate computations,
one does \emph{not} compute $\log(p^k)=k\log(p)$, but only some rough
approximation, for instance by counting the number of bytes that the
multi-precision integer $p^k$ occupies in memory, or any other fast method.

Even though this method is designed to be fast for relatively large $k$,
we find that it is considerably faster than any of the Bernoulli methods, even
for very small $k$, the ratio increasing with increasing $k$ and decreasing 
$F$.

Here are the times obtained using this method. The reader will notice that the
times for very small $k$ are larger than for moderate $k$ due to the very
large number of Euler factors to be computed, the smallest being impossibly
long. We use $*$ to indicate very long times (usually more than $100$ seconds),
and on the contrary -- to indicate a negligible time, less than $50$
milliseconds.

\bigskip

\centerline{
\begin{tabular}{|r||r|r|r|r|r|r|r|r|r|}
\hline
$D\dd k$ & $2$    & $4$    & $6$    & $8$    & $10$   & $12$   & $14$   & $16$ & $18$ \\
\hline\hline
$10^6+1$ & $*$    & $1.03$ & $0.20$ & $0.11$ & $0.09$ & $0.08$ & $0.08$ & $0.07$ & $0.08$  \\
$10^7+1$ & $*$    & $14.4$ & $2.36$ & $1.19$ & $0.93$ & $0.81$ & $0.79$ & $0.73$ & $0.77$ \\
$10^8+1$ & $*$    & $*$    & $27.9$ & $13.1$ & $9.75$ & $8.32$ & $7.97$ & $7.25$ & $7.54$ \\
$10^9+1$ & $*$    & $*$    & $*$    & $*$    & $105.$ & $87.2$ & $81.7$ & $73.5$ & $75.5$ \\
\hline
\end{tabular}}

\bigskip

\centerline{
\begin{tabular}{|r||r|r|r|r|r|r|r|r|r|}
\hline
$D\dd k$ & $20$   & $40$   & $60$   & $80$   & $100$   & $150$   & $200$   & $250$ & $300$ \\
\hline\hline
$10^5+1$ &  --    &  --    &  --    &  --    &  --     & $0.06$  & $0.08$  & $0.11$ & $0.15$  \\
$10^6+1$ & $0.08$ & $0.12$ & $0.17$ & $0.22$ & $0.29$ & $0.48$ & $0.68$ & $1.01$ & $1.32$  \\
$10^7+1$ & $0.77$ & $1.09$ & $1.62$ & $2.01$ & $2.66$ & $4.48$ & $6.29$ & $9.23$ & $12.2$ \\
$10^8+1$ & $7.52$ & $10.3$ & $15.1$ & $18.8$ & $24.7$ & $41.3$ & $58.5$ & $85.5$ & $114.$ \\
$10^9+1$ & $75.2$ & $99.8$ & $*$ & $*$ & $*$ & $*$ & $*$ & $*$ & $*$ \\
\hline
\end{tabular}}

\bigskip

\bigskip

\centerline{
\begin{tabular}{|r||r|r|r|r|r|r|r|r|}
\hline
$D\dd k$ & $400$  & $800$  & $1600$ & $3200$ & $6400$ & $12800$ & $25600$ & $51200$ \\
\hline\hline
$10^2+5$ & --     & --     & --     &  --    & $0.24$ & $1.25$ & $6.36$ & $31.3$ \\
$10^3+1$ & --     & --     & $0.07$ & $0.34$ & $1.85$ & $9.84$ & $51.2$ & $*$  \\
$10^4+1$ & --     & $0.11$ & $0.56$ & $2.94$ & $15.9$ & $85.6$ & $*$    & $*$ \\
$10^5+1$ & $0.24$ & $0.98$ & $4.96$ & $26.1$ & $*$    & $*$    & $*$    & $*$  \\
$10^6+1$ & $2.16$ & $9.49$ & $44.4$ & $*$    & $*$    & $*$    & $*$    & $*$\\
$10^7+1$ & $20.0$ & $82.3$ & $*$    & $*$    & $*$    & $*$    & $*$    & $*$\\
\hline
\end{tabular}}

\bigskip

\section{Using the Approximate Functional Equation}

All the methods that we have seen up to now take time proportional to the
conductor $F$, the main difference between them being the dependence in $k$
and the size of the implicit constant in the time estimates.

We are now going to study a number of methods which take time proportional
to $F^{1/2+\eps}$ for any $\eps>0$. The simplest version of the
\emph{approximate functional equation} that we will use is as follows:

\begin{theorem} Let $e=0$ or $1$ be such that $\chi(-1)=(-1)^k=(-1)^e$.
For any complex $s$ we have the following formula, valid for any $A>0$:
$$\G\left(\dfrac{s+e}{2}\right)L(\chi,s)=
\sum_{n\ge1}\dfrac{\chi(n)}{n^s}\G\left(\dfrac{s+e}{2},\dfrac{A\pi n^2}{F}\right)+\eps(\chi)\sum_{n\ge1}\dfrac{\ov{\chi}(n)}{n^{1-s}}\G\left(\dfrac{1-s+e}{2},\dfrac{\pi n^2}{AF}\right)\;,$$
where the \emph{root number} $\eps(\chi)$ is given by the formula
$\eps(\chi)=\g(\chi)/(i^e\sqrt{F})$, where $\g(\chi)$ is the Gauss sum
attached to $\chi$, and $\G(s,x)$ is the incomplete gamma function
$\G(s,x)=\int_x^\infty t^{s-1}e^{-t}\,dt$.
\end{theorem}

Since $\G(s,x)\sim x^{s-1}e^{-x}$ hence tends to $0$ exponentially
fast as $x\to\infty$, the above formula does lead to a $\Os(F^{1/2})$
algorithm for computing $L(\chi,s)$, not necessarily for a negative
integer $s$. Note that this type of formula is available for any
type of $L$-function with functional equation, not only those attached
to a Dirichlet character.

The constant $A$ is included as a check on the implementation, since the
left-hand side is independent of $A$, but once checked the optimal choice
is $A=1$. This constant can also be used to compute $\eps(\chi)$ if it
is not known (note that $\eps(\chi)=1$ if $\chi$ is quadratic), but there
are better methods to do this.

Even though this method is in $\Os(F^{1/2})$, so asymptotically much faster
than the first two methods that we have seen, its main drawback is the
computation time of $\G(s,x)$. Even though quite efficient methods are
known for computing it, our timings have shown that in all ranges of
the conductor $F$ and value of $k$, either the use of the full functional
equation or the use of Eisenstein series of half-integral weight (methods
(2) and (5)) are considerably faster, so we will not discuss this method
further.

\section{Using Hecke--Eisenstein Series}

\subsection{The Main Theorem}

The main theorem comes from the computation of the Fourier expansion of
Hecke--Eisenstein series in the theory of Hilbert modular forms, and
is easily proved using the methods of \cite{CohH}:

\begin{theorem}\label{mainth} Let $K$ be a real quadratic field of 
discriminant $D>0$, let $\psi$ be a primitive character modulo $F$ such that 
$\psi(-1)=(-1)^k$, let $N$ be a squarefree integer, and assume that 
$\gcd(F,ND)=1$. If we set
$$a_{k,\psi,N}(0)=\prod_{p\mid N}(1-\psi\chi_D(p)p^{k-1})\dfrac{L(\psi,1-k)L(\psi\chi_D,1-k)}{4}\;,$$
and for $n\ge1$:
$$a_{k,\psi,N}(n)=\sum_{\substack{d\mid n\\\gcd(d,N)=1}}\psi\chi_D(d)d^{k-1}\sum_{s\in\Z}\sigma_{k-1,\psi}\left(\dfrac{(n/d)^2D-s^2}{4N}\right)\;,$$
where $\sigma_{k-1,\psi}(m)=\sum_{d\mid m}\psi(d)d^{k-1}$, then
$$\sum_{n\ge0}a_{k,\psi,N}(n)q^n\in M_{2k}(\G_0(FN),\psi^2)\;.$$
\end{theorem}

Note that in the above we set implicitly $\sigma_{k-1,\psi}(m)=0$ if
$m\notin\Z_{\ge1}$.

The restriction $\gcd(F,N)=1$ is not important, since letting $N$ have factors
in common with $F$ would not give more general results. Similarly for the
restriction on $N$ being squarefree. On the other hand,
the restriction $\gcd(F,D)=1$ is more important: similar results exist
when $\gcd(F,D)>1$, but they are considerably more complicated. Since we need
them, we will give one such result below in the case $\gcd(F,D)=4$. 

\smallskip

We use this theorem in the following way. First, we must assume for
practical reasons that $k$, $F$, and $N$ are not too large. In that case it is
very easy to compute explicitly a basis for $M_{2k}(\G_0(FN),\psi^2)$.
Given this basis, it is then easy to express any constant term of an element
of the space as a linear combination of $u$ coefficients (not necessarily the
first ones), where $u$ is the dimension of the space. In particular, this
gives $a_{k,\psi,N}(0)$, and hence $L(\psi\chi_D,1-k)$, as a finite linear
combination of some $a_{k,\psi,N}(n)$ for $n\ge1$.

Second, since the conductor of $\psi$ must be small, the method is thus
applicable only to compute $L(\chi,1-k)$ for Dirichlet characters $\chi$
which are ``close'' to quadratic characters, in other words of the form
$\psi\chi_D$ with conductor of $\psi$ small.
Of course the quantities $L(\psi,1-k)$ are computed once and for all (using
any method, since $F$ and $k$ are small). Note that the auxiliary
integer $N$ is used only to improve the speed of the formulas, as we will
see below, but of course one can always choose $N=1$ if desired.

\subsection{The Case $k$ Even}

For future use, define
$$S_k(m,N)=\sum_{s\in\Z}\sigma_k\left(\dfrac{m-s^2}{N}\right)\;,$$
where for any arithmetic function $f$ such as $\sigma_k$ we set
$f(x)=0$ if $x\not\in\Z_{\ge1}$, i.e., if $x$ is either not integral or
non-positive. Using the theorem with $F=N=1$ we immediately obtain formulas
such as
\begin{align*}
  L(\chi_D,-1)&=-\dfrac{1}{5}S_1(D,4)\\
  L(\chi_D,-3)&=S_3(D,4)\\
  L(\chi_D,-5)&=-\dfrac{1}{195}\left(\left(24+2^5\leg{D}{2}\right)S_5(D,4)+S_5(D,1)\right)\\
\end{align*}

To obtain a general formula we recall the following theorem of Siegel:
\begin{theorem}\label{thm:siegel} Let $r=\dim(M_k(\G))$ and define coefficients
$c_i^k$ by
$$\Delta^{-r}E_{12r-k+2}=\sum_{i\ge-r}c_i^kq^i\;,$$
where by convention $E_0=1$. Then for any $f=\sum_{n\ge0}a(n)q^n\in M_k(\G)$
we have
$$\sum_{0\le n\le r}a(n)c_{-n}^k=0\;,$$
and $c_0^k\ne0$.\end{theorem}

Combined with the main theorem (with $F=N=1$), we obtain the following
corollary:

\begin{corollary}\label{cor:siegel} Keep the above notation, let $k\ge2$ be an
even integer, and set $r=\dim(M_{2k}(\G))=\lfloor k/6\rfloor+1$. If $D>0$ is a
fundamental discriminant we have
$$L(\chi_D,1-k)=\dfrac{4k}{c_0^{2k}B_k}\sum_{1\le m\le r}S_{k-1}(m^2D,4)\sum_{1\le d\le r/m}d^{k-1}\leg{D}{d}c_{-dm}^{2k}\;.$$
\end{corollary}

For very small values of $k$ it is possible to improve on the speed of
the above general formula by choosing $F=1$ but larger values of $N$ in the
theorem. Without entering into details, on average we can gain a factor of
$3.95$ for $k=2$, of $1.6$ for $k=6$, and of $1.1$ for $k=8$, and I have
found essentially no improvement for other values of $k$ including for $k=4$.

The advantages of this method are threefold. First, it is by far the fastest
method seen up to now for computing $L(\chi_D,1-k)$. Second, the universal
coefficients $c_{-n}^k$ that we need are easily computed thanks to Siegel's
theorem. And third, the flexibility of choosing the auxiliary Dirichlet
character $\psi$ in the theorem allows us to compute $L(\chi,1-k)$ for
more general characters $\chi$ than quadratic ones.

The two disadvantages are first that the quantities $S_{k-1}(m^2D,4)$ need to
be computed for each $m$ (although some duplication can be avoided),
and second that $m^2D$ becomes large when $m$ increases. These two
disadvantages will disappear in the method using Eisenstein series of
half-integral weight (at the expense of losing some of the advantages
mentioned above), so we will not give the timings for this method.

\subsection{The Case $k$ Odd}

Thanks to the main theorem, although Hilbert modular forms in two variables are
only for \emph{real} quadratic fields, thus with discriminant $D>0$, if we
choose an odd character $\psi$ such as $\psi=\chi_{-3}$ or $\chi_{-4}$, it can
also be used to compute $L(\chi_D,1-k)$ for $D<0$, hence $k$ odd. I have not
been able to find useful formulas with $\psi=\chi_{-3}$, so from now on we
assume that $\psi=\chi_{-4}$, so $F=4$. We first introduce some notation.

\begin{definition}\label{defsigchi}
We set
$$\sigma_k^{(1)}(m)=\sum_{d\mid m}\leg{-4}{d}d^k\;,\quad\sigma_k^{(2)}(m)=\sum_{d\mid m}\leg{-4}{m/d}d^k\;,\text{\quad and}$$
$$S_k^{(j)}(m,N)=\sum_{s\in\Z}\sigma_k^{(j)}\left(\dfrac{m-s^2}{N}\right)\;,$$
with the usual understanding that $\sigma_k^{(j)}(m)=0$ if $m\notin\Z_{\ge1}$.
\end{definition}

Note that, as for $\sigma_k$ itself when $k$ is odd, for $k$ even these
arithmetic functions occur naturally as Fourier coefficients of Eisenstein
series of weight $k+1$ and character $\lgs{-4}{.}$. More precisely, for
$k\ge3$ odd the series $E_k(\chi_{-4},1)$ and $E_k(1,\chi_{-4})$ form a basis
of the Eisenstein subspace of $M_k(\G_0(4),\chi_{-4})$, where
\begin{align*}
E_k(\chi_{-4},1)(\tau)&=\dfrac{L(\chi_{-4},1-k)}{2}+\sum_{n\ge1}\sigma_{k-1}^{(1)}(n)q^n\text{\quad and}\\
E_k(1,\chi_{-4})(\tau)&=\sum_{n\ge1}\sigma_{k-1}^{(2)}(n)q^n\;.\end{align*}

\smallskip

To be able to use the theorem in general, it is necessary to assume the
following:

\begin{conjecture} If $D>1$ is squarefree (not necessarily a discriminant),
$F=4$, and $N=1$, the statement of Theorem \ref{mainth} is still valid
verbatim.\end{conjecture}

This is probably easy to prove, and I have checked it on thousands of examples.
Assuming this conjecture, applying the theorem to $\psi=\chi_{-4}$ and
the Hecke operator $T(2)$ it is immediate to prove the following:

\begin{corollary}\label{cor2}
 Let $D<-4$ be any fundamental discriminant. Set
$$a_{k,D}(0)=\left(1-2^{k-1}\leg{D}{2}\right)\dfrac{L(\chi_{-4},1-k)L(\chi_D,1-k)}{4}\;,\text{\quad and}$$
 $$a_{k,D}(n)=\sum_{d\mid n}\leg{4D/\delta}{d}d^{k-1}S_{k-1}^{(1)}((n/d)^2|D/\delta|,1)\;,$$
where $\delta=1$ if $D\equiv1\pmod4$ and $\delta=4$ if $D\equiv0\pmod4$.
Then $\sum_{n\ge0}a_{k,D}(n)q^n\in M_{2k}(\G_0(2))$.\end{corollary}

To use this result, we need an analogue of Siegel's Theorem \ref{thm:siegel}
for $\G_0(2)$, and for this we need to introduce a number of modular forms.
\begin{definition}
We set $F_2(\tau)=2E_2(2\tau)-E_2(\tau)$, $F_4(\tau)=(16E_4(2\tau)-E_4(\tau))/15$,
and $\D_4(\tau)=(E_4(\tau)-E_4(2\tau))/240$, where $E_2$ and $E_4$ are the
standard Eisenstein series of weight $2$ and $4$ on the full modular group.
\end{definition}
Note that $F_2\in M_2(\G_0(2))$ and $F_4$ and $\D_4$ are in $M_4(\G_0(2))$.

\begin{theorem} Let $k\in2\Z$ be a positive even integer, set
$r=\lfloor k/4\rfloor+2$, $E=F_2F_4$ if $k\equiv0\pmod4$, $E=F_4$ if
$k\equiv2\pmod4$, and write $E/\D_4^r=\sum_{i\ge-r}c_i^kq^i$. Then for any 
$F=\sum_{n\ge0}a(n)q^n\in M_k(\G_0(2))$ we have
$$\sum_{0\le n\le r}a(n)c_{-n}^k=0\;,$$
and in addition $c_0^k\ne0$.
\end{theorem}

Note that since we will use this theorem for $M_{2k}(\G_0(2))$ with $k$ odd,
we have $2k\equiv2\pmod4$, so we will always use $E=F_4$. The analogue of
Corollary \ref{cor:siegel} is then as follows:

\begin{corollary} Keep the above notation, let $k\ge3$ be an odd integer, and
set $r=(k+3)/2$. If $D<-4$ is a fundamental discriminant we have
$$L(\chi_D,1-k)=\dfrac{8}{A}\sum_{1\le m\le r}S^{(1)}_{k-1}(m^2|D|/\delta,1)\sum_{1\le d\le r/m}d^{k-1}\leg{4D/\delta}{d}c_{-dm}^{2k}\;,$$
with $A=c_0^{2k}(2^{k-1}\lgs{D}{2}-1)E_{k-1}$, and
where the $E_k$ are the \emph{Euler numbers} ($E_0=1$, $E_2=-1$, $E_4=5$,
$E_6=-61$,\dots).
\end{corollary}

The advantages/disadvantages mentioned in the case $k$ even are the same
here.

\section{Using Eisenstein Series of Half-Integral Weight}

We now come to the most powerful method known to the author for computing
$L(\chi_D,1-k)$: the use of Eisenstein series of half-integral weight. Once
again, we will see a sharp distinction between $k$ even and $k$ odd. We
first begin by recalling some basic results on $M_w(\G_0(4))$ (we use
the index $w$ for the weight since it will be used with $w=k+1/2$). Later, we
will see that it is more efficient to use modular forms on subgroups of
$\G_0(4)$.
  
\subsection{Basic Results on $M_w(\G_0(4))$}

Recall that the basic theta function
$\th(\tau)=\sum_{s\in\Z}q^{s^2}=1+2\sum_{s\ge1}q^{s^2}$
satisfies for any $\ga=\psmm{a}{b}{c}{d}\in\G_0(4)$ the modularity condition
$\th(\ga(\tau))=v_{\th}(\ga)(c\tau+d)^{1/2}\th(\tau)$,
where the \emph{theta-multiplier system} $v_{\th}(\ga)$ is given by
$$v_{\th}(\ga)=\leg{-4}{d}^{-1/2}\leg{c}{d}\;,$$
and all square roots are taken with the principal determination.
The space $M_w(\G_0(4),v_{\th}^{2w})$ of holomorphic functions behaving
modularly like $\th^{2w}$ under $\G_0(4)$ and holomorphic at the cusps will be
simply denoted $M_w(\G_0(4))$ since there is no risk of confusion. Note,
however, that if $w$ is an odd integer and in the context of modular forms
of integral weight, $M_w(\G_0(4))$ is denoted $M_w(\G_0(4),\chi_{-4})$.

We recall the following easy and well-known results (note that $F_2$ and
$\D_4$ are not the same functions as those used above):

\begin{proposition}
Define
\begin{align*}
F_2(\tau)&=\dfrac{\eta(4\tau)^8}{\eta(2\tau)^4}=-\dfrac{1}{24}(E_2(\tau)-3E_2(2\tau)+2E_2(4\tau))\;,\\
\D_4(\tau)&=\dfrac{\eta(\tau)^8\eta(4\tau)^8}{\eta(2\tau)^8}=\dfrac{1}{240}(E_4(\tau)-17E_4(2\tau)+16E_4(4\tau))\;.\end{align*}
\begin{enumerate}
\item We have
$$\bigoplus_wM_w(\G_0(4))=\C[\th,F_2]\text{\quad and\quad}
\bigoplus_wS_w(\G_0(4))=\th\D_4\C[\th,F_2]\;.$$
\item In particular we have the dimension formulas
$$\dim(M_w(\G_0(4)))=\begin{cases}0&\text{\quad for $w<0$}\\
1+\lfloor w/2\rfloor&\text{\quad for $w\ge0$\;.}\end{cases}$$
$$\dim(S_w(\G_0(4)))=\begin{cases}0&\text{\quad for $w\le4$}\\
\lfloor w/2\rfloor-1&\text{\quad for $w>2$, $w\notin2\Z$}\\
\lfloor w/2\rfloor-2&\text{\quad for $w>2$, $w\in2\Z$\;.}
\end{cases}$$\end{enumerate}
\end{proposition}

We also recall that when $w\in1/2+\Z$, the Kohnen $+$-space of $M_w(\G_0(4))$,
denoted simply by $M_w^+$, is defined to be the space of forms $F$ having a
Fourier expansion $F(\tau)=\sum_{n\ge0}a(n)q^n$ with $a(n)=0$ if
$(-1)^{w-1/2}n\not\equiv0,1\pmod4$. Note that we include Eisenstein series.
It is clear that $M_{1/2}^+=M_{1/2}(\G_0(4))$ and $M_{3/2}^+=\{0\}$, so
we will always assume that $w\ge5/2$. In that case there a single
Eisenstein series in $M_w^+$, due to the author, that we will denote by
$\H_k$: its importance is due to the fact that if we write
$\H_k(\tau)=\sum_{n\ge0}a_k(n)q^n$, then if $D=(-1)^{w-1/2}n$ is a fundamental
discriminant we have $a_k(n)=L(\chi_D,1-(w-1/2))$, so being able to compute
efficiently the Fourier coefficients of $\H_k$ automatically gives us a fast
method for computing our desired quantities $L(\chi_D,1-k)$ with $k=w-1/2$.

The remaining part of $M_w^+$, which we of course denote by $S_w^+$,
is formed by the cusp forms belonging to $M_w^+$. One of Kohnen's main theorems
is that $S_w^+$ is Hecke-isomorphic to the space of modular forms of even
weight $S_{2w-1}(\G)$. In particular, note the following:

\begin{corollary} For $w\ge5/2$ half-integral we have
$$\dim(M_w^+)=\begin{cases}
1+\lfloor w/6\rfloor&\text{\quad if $6\nmid(w-3/2)$\;,}\\
\lfloor w/6\rfloor&\text{\quad if $6\mid(w-3/2)$\;.}\end{cases}$$
\end{corollary}

\smallskip

Notation: \begin{enumerate}
\item Recall that if $a(n)$ is any arithmetic function (typically
$a=\sigma_{k-1}$ or twisted variants), we define $a(x)=0$ if
$x\notin\Z_{\ge1}$.
\item If $F$ is a modular form and $d\in\Z_{\ge1}$, we denote by $F[d]$
the function $F(d\tau)$.
\item We will denote by $D(F)$ the differential operator 
$qd/dq=(1/(2\pi i))d/d\tau$.\end{enumerate}

\subsection{The Case $k$ Even using $\G_0(4)$}

The main idea is to use Rankin--Cohen brackets of known series with $\th$ to
generate $M_w^+$: indeed, $\th$ and its derivatives
are \emph{lacunary}, so multiplication by them is much faster than ordinary
multiplication, at least in reasonable ranges (otherwise Karatsuba or FFT
type methods are faster to construct \emph{tables}).

First note the following immediate result:

\begin{proposition} The form $\th E_2[4]-6D(\th)$ is a basis of $M_{5/2}^+$
and the form $\th E_4[4]$ is a basis of $M_{9/2}^+$.\end{proposition}

In particular, we recover the formulas $L(\chi_D,-1)=(-1/5)S_1(D,4)$
and $L(\chi_D,-3)=S_3(D,4)$ already obtained using Hecke--Eisenstein series.

\smallskip

As in the case of Hecke--Eisenstein series, we will need to distinguish two
completely different cases: the case $w-1/2$ \emph{even}, which is considerably
simpler, and the case $w-1/2$ odd, which is more complicated and less
efficient. The reason for this sharp distinction is the following theorem:

\begin{theorem} Assume that $w\ge9/2$ is such that $k=w-1/2\in2\Z$. The
modular forms 
$$([\th,E_{k-2j}[4]]_j)_{0\le j\le\lfloor k/6\rfloor}$$
form a basis of $M_w^+$, where we recall that $[f,g]_n$ denotes the
$n$th Rankin--Cohen bracket.\end{theorem}

We can now easily achieve our goal. First, we compute the Fourier
expansions of the basis given by the theorem up to the Sturm bound.
Then to compute $L(\chi_D,1-k)$ with $k=w-1/2$, we do as follows: we compute
the Fourier expansion of $\H_k$ up to the Sturm bound, and using the basis
coefficients we deduce a linear combination of the form
$$\H_k=\sum_{0\le j\le \lfloor k/6\rfloor}c_j^k[\th,E_{k-2j}[4]]_j\;.$$
We can easily compute the coefficients of these brackets:

\begin{proposition}\label{prop:even4} Let $F_r=-B_rE_r/(2r)$ be the Eisenstein
  series of level
  $1$ and weight $r$ normalized so that the Fourier coefficient $q^1$ is
  equal to $1$. We have $[\th,F_r[4]]_n=\sum_{m\ge0}b_{n,r}(m)q^m\;,$ with
  \begin{align*}b_{n,r}(m)&=m^n\sum_{s\in\Z}P_{n,r}(s^2/m)\sigma_{r-1}\left(\dfrac{m-s^2}{4}\right)\text{ , where}\\
    P_{n,r}(X)&=\sum_{\ell=0}^n(-1)^\ell\binom{n-1/2}{\ell}\binom{2n+r-\ell-3/2}{n-\ell}X^{n-\ell}\;,\end{align*}
are \emph{Gegenbauer polynomials}.\end{proposition}

In particular, if we generalize a previous notation and set for any
polynomial $P_n$ of degree $n$
$$S_k(m,N,P_n)=m^n\sum_{s\in\Z}P_n(s^2/m)\sigma_k\left(\dfrac{m-s^2}{N}\right)\;,$$
we have
$$L(\chi_D,1-k)=\sum_{0\le j\le\lfloor k/6\rfloor}c_j^kS_{k-2j-1}(D,4,P_{j,k-2j})\;.$$
The biggest advantage of this formula compared to the one coming from
Hecke--Eisenstein series is that the different $S_{k-2j-1}$ can be computed
simultaneously since they involve factoring the same integers $(D-s^2)/4$,
and in addition these integers stay small, contrary to the former method
where the integers were of the form $(n^2D-s^2)/4$.

The two disadvantages are that first, it is not easy (although possible)
to generalize to general characters $\chi$, but mainly because for large $k$
the computation of $c_j^k$ involves solving a linear system of
size proportional to $k$, so when $k$ is in the thousands, this becomes
prohibitive. It is possible that there is a much faster method to compute
them analogous to Siegel's theorem which expresses the constant term of a
modular form as a universal (for a given weight) linear combination of higher
degree terms, but I do not know of such a method.

As already mentioned, this gives the fastest method known
to the author for computing $L(\chi_D,1-k)$, at least when $k$ is not
unreasonably large.

\subsection{The Case $k$ Even using $\G_0(4N)$}

We can, however, do better by using subgroups of $\G_0(4)$, i.e., brackets
with $E_{k-2j}[4N]$ for $N>1$. Recall that in the case of Hecke--Eisenstein
series this allowed us to give faster formulas only for very small values
of $k$ ($k=2$, $6$ and $8$). On the contrary, we are going to see that
here we can obtain faster formulas for all $k$, only depending on
congruence and divisibility properties of the discriminant $D$.

After considerable experimenting, I have arrived at the following
conjecture, which I have tested on tens of thousands of cases and
\emph{proved} in small weights. All of these identities can in principle
be proved.

\begin{conjecture}\label{conjeven} For $N=4$, $8$, $12$ and $16$ and any even
integer $k\ge2$ there exist universal coefficients $c_{j,N}^k$ such that for
all positive fundamental discriminants $D$ (which in addition must be
congruent to $1$ modulo $8$ when $N=16$) we have
$$\left(1+\leg{D}{N/4}\right)L(\chi_D,1-k)=\sum_{0\le j\le\lfloor k/m_N\rfloor}c_{j,N}^kS_{k-2j-1}(D,N,P_{j,k-2j})\;,$$
with $m_N=6$, $4$, $3$, and $4$ respectively and the same polynomials
$P$ as above.\end{conjecture}

By what we said above this conjecture is proved for $N=4$ (with
$c_{j,4}^k=2c_j^k$), and should be easy to prove using the
finite-dimensionality of the corresponding modular form spaces together with
the Sturm bounds, but I have not done these proofs. It is also easy to prove
for small values of $k$.

It is clear that if we can choose a larger value of $N$ than $N=4$ (i.e.,
when $1+\lgs{D}{N/4}\ne0$) the number of terms involved in $S_{k-2j-1}$
will be smaller. Computing that number leads to the following algorithm:

If $3\mid D$ use $N=12$, otherwise if $D\equiv1\pmod8$ use $N=16$, otherwise
if $4\mid D$ use $N=8$, otherwise if $D\equiv1\pmod3$ use $N=12$, otherwise
use $N=4$.

Note, however, that the size of the linear system which needs to be solved
to find the coefficients $c_{j,N}^k$ is larger when $N>4$, so one must
balance the time to compute these coefficients with the size of $D$: for very
large $D$ it may be worthwhile, but for moderately large $D$ it may be
better to always choose $N=4$ (see the second table below).

\medskip

As before, we give tables of timings using these improvements. Note that
they are only an indication, since congruences modulo $16$ and $3$
may improve the times:

\bigskip

\centerline{
\begin{tabular}{|r||r|r|r|r|r|r|r|r|r|}
\hline
$D\dd k$ & $2$    & $4$    & $6$    & $8$    & $10$   & $12$   & $14$   & $16$ & $18$ \\
\hline\hline
$10^{10}+1$ & $0.07$ & $0.07$ & $0.07$ & $0.08$ & $0.08$ & $0.09$ & $0.09$ & $0.11$ & $0.11$ \\
$10^{11}+9$ & $0.30$ & $0.32$ & $0.33$ & $0.35$ & $0.36$ & $0.39$ & $0.40$ & $0.44$ & $0.44$ \\
$10^{12}+1$ & $2.25$ & $2.31$ & $2.32$ & $2.41$ & $2.42$ & $2.53$ & $2.55$ & $2.67$ & $2.69$ \\
$10^{13}+1$ & $10.3$ & $10.5$ & $10.5$ & $10.8$ & $10.9$ & $11.2$ & $11.3$ & $11.7$ & $11.8$ \\
$10^{14}+1$ & $54.0$ & $54.7$ & $55.0$ & $55.8$ & $56.2$ & $57.3$ & $57.6$ & $59.0$ & $59.3$ \\
\hline
\end{tabular}}

\bigskip

In the next table, we use the improvements for larger $N$ only when $D$ is
sufficiently large, and the corresponding timings have a ${}^*$; all the others
are obtained only with $N=4$:

\bigskip

\centerline{
\begin{tabular}{|r||r|r|r|r|r|r|r|r|r|r|}
\hline
$D\dd k$ & $20$   & $40$   & $60$   & $80$   & $100$  & $150$  & $200$  & $250$  & $300$ & $350$ \\
\hline\hline
$10^6+1$ & --     & --     & --     & --     & --     & --     & $0.07$ & $0.14$ & $0.29$ & $0.51$ \\
$10^7+1$ & --     & --     & --     & --     & --     & $0.08$ & $0.16$ & $0.30$ & $0.55$ & $0.88$ \\
$10^8+1$ & --     & --     & --     & $0.06^*$ & $0.10^*$ & $0.25$ & $0.50$ & $0.89$ & $1.50$ & $2.28$ \\
$10^9+1$ & --     & $0.06^*$ & $0.12^*$ & $0.19^*$ & $0.30^*$ & $0.74^*$ & $1.59^*$ & $2.85^*$ & $4.96^*$ & $7.56$ \\
$10^{10}+1$ & $0.12^*$ & $0.23^*$ & $0.39^*$ & $0.64^*$ & $1.00^*$ & $2.48^*$ & $5.20^*$ & $9.08^*$ & $15.0^*$ & $22.4^*$ \\
$10^{11}+9$ & $0.49^*$ & $0.86^*$ & $1.47^*$ & $2.37^*$ & $3.67^*$ & $9.09^*$ & $18.9^*$ & $32.8^*$ & $52.8^*$ & $77.2^*$ \\
$10^{12}+1$ & $2.84^*$ & $4.04^*$ & $6.01^*$ & $9.04^*$ & $13.4^*$ & $31.8^*$ & $64.6^*$ & $*$ & $*$ & $*$ \\
$10^{13}+1$ & $12.3^*$ & $16.5^*$ & $23.4^*$ & $34.2^*$ & $49.9^*$ & $*$ & $*$ & $*$ & $*$ & $*$ \\
$10^{14}+1$ & $60.8^*$ & $74.8^*$ & $98.8^*$ & $*$ & $*$ & $*$ & $*$ & $*$ & $*$ \\
\hline
\end{tabular}}

\bigskip

For larger values of $k$ the time to compute the coefficients dominate,
so we first give a table giving these timings:

\bigskip

\centerline{
\begin{tabular}{|r||r|r|r|r|r|r|r|r|r|r|}
\hline
$N\dd k$ & $100$ & $200$ & $300$ & $400$   & $500$   & $600$   & $700$   & $800$   & $900$ & $1000$ \\
\hline\hline
   $4$ & -- & $0.04$ & $0.20$ & $0.69$ & $1.95$ & $4.04$ & $7.54$ & $13.3$ & $22.4$ & $34.4$ \\
   $8$ & -- & $0.17$ & $0.87$ & $2.77$ & $6.95$ & $14.7$ & $28.4$ & $49.2$ & $83.5$ & $*$ \\
  $12$ & -- & $0.32$ & $1.90$ & $5.77$ & $14.5$ & $32.0$ & $61.6$ & $*$ & $*$ & $*$ \\
  $16$ & -- & $0.20$ & $1.13$ & $3.64$ & $9.59$ & $20.5$ & $31.4$ & $53.4$ & $89.8$ & $*$ \\
\hline
\end{tabular}}

\bigskip
  
As already mentioned, these timings would become much smaller if
we had a method analogous to Siegel's theorem to compute them.

\bigskip

\centerline{
\begin{tabular}{|r||r|r|r|r|r|r|r|r|}
\hline
$D\dd k$ & $400$   & $500$   & $600$   & $700$   & $800$   & $900$ & $1000$ \\
\hline\hline
$10^5+1$ & $0.73$  & $2.03$  & $4.16$  & $7.72$  & $13.6$  & $22.7$ & $34.8$ \\
$10^6+1$ & $0.87$  & $2.26$  & $4.53$  & $8.27$  & $14.3$  & $23.8$ & $36.2$ \\
$10^7+1$ & $1.39$  & $3.21$  & $6.91$  & $10.4$  & $17.4$  & $27.8$ & $41.6$ \\
$10^8+1$ & $3.33$  & $6.61$  & $11.4$  & $18.3$  & $28.4$  & $42,7$ & $61.7$ \\
$10^9+1$ & $10.7$  & $19.6$  & $32.0$  & $48.1$  & $70.1$  & $99.3$ & $*$    \\
$10^{10}+1$ & $31.9^*$  & $58.9^*$  & $98.9^*$  & $*$     & $*$     & $*$ & $*$  \\
$10^{11}+9$ & $108.^*$ & $*$  & $*$  & $*$  & $*$  & $*$  & $*$ \\
\hline
\end{tabular}}

\bigskip

\subsection{The Case $k$ Odd using $\G_0(4N)$}

In this case, the Kohnen $+$-space, to which $\H_k$ belongs, is the
space of modular forms $\sum_{n\ge0}a(n)q^n$ such that $a(n)=0$ if
$n\equiv1$ or $2$ modulo $4$. Thus, we cannot hope to \emph{directly} obtain
elements in this space using brackets with $\th$. What we can do is the
following: as above, for $\ell\ge1$ odd consider the two Eisenstein series
\begin{align*}
E_{\ell}^{(1)}:=E_{\ell}(\chi_{-4},1)(\tau)&=\dfrac{L(\chi_{-4},1-\ell)}{2}+\sum_{n\ge1}\sigma_{\ell-1}^{(1)}(n)q^n\text{\quad and}\\
E_{\ell}^{(2)}:=E_{\ell}(1,\chi_{-4})(\tau)&=\sum_{n\ge1}\sigma_{\ell-1}^{(2)}(n)q^n\;,\end{align*}
which belong to $M_{\ell}(\G_0(4))$ (using our notation, otherwise we should
write $M_{\ell}(\G_0(4),\chi_{-4})$). It is clear that for $u=1$ and $2$ the
$j$-th brackets $[\th,E_{k-2j}^{(u)}]_j$ belong to $M_{k+1/2}(\G_0(4))$, and
it should be easy to prove that they generate this space (I have extensively
tested this, and if it was not the case the implementation would detect it).
We can therefore express any modular form, in particular $\H_k$, as
a linear combination of these brackets, and therefore again obtain explicit
formulas for $L(\chi_D,1-k)$.

However, we can immediately do considerably better in two different ways.
First, by Shimura theory we know that $T(4)\H_k$ still belongs to
$M_{k+1/2}(\G_0(4))$, and by definition it is equal to
$\sum_{n\ge0}H_k(4n)q^n$. Expressing it as a linear combination of the
above brackets again gives formulas for $L(\chi_D,1-k)$, but where the
coefficients involve $|D|/4$ instead of $|D|$, so much faster (and of course
applicable only for $D\equiv0\pmod4$). Note that this trick is \emph{not}
applicable in the case of even $k$ because $T(4)\H_k$ is not anymore
in the Kohnen $+$-space, so we would lose all the advantages of having a space
of small dimension.

The second way in which we can do better is to consider brackets of $\th$ with
$E_{\ell}^{(u)}[N]$ (where we replace $q^n$ by $q^{Nn}$) for suitable values of
$N$: note that these modular forms belong to $M_{k+1/2}(\G_0(4N))$. It is then
necessary to apply a Hecke-type operator to reduce the dimension of the
spaces that we consider. More precisely, if we only look at coefficients
$a(|D|)$ with given $\lgs{D}{2}$, we see experimentally that there is a
linear relation between $\H_k$ and the above brackets. This leads to the
following analogue for $k$ odd of Conjecture \ref{conjeven}, where
generalizing the notation $S_k^{(j)}(m,N)$ used above for $j=1$ and $2$ we
also use
$$S_k^{(j)}(m,N,P_n)=m^n\sum_{s\in\Z}P_n(s^2/m)\sigma_k^{(j)}\left(\dfrac{m-s^2}{N}\right)\;,$$
where $P_n$ is a polynomial of degree $n$.

\begin{conjecture}\label{conjodd} For $N=1$, $2$, $3$, $5$, $6$, and $7$,
any odd integer $k\ge3$, and $e\in\{-1,0,1\}$, there exist universal
coefficients $c_{j,N,e}^k$ such that for all negative fundamental discriminants
$D$ such that $\lgs{D}{2}=e$ we have
$$\left(1+\leg{|D|}{N_2}\right)L(\chi_D,1-k)=\sum_{0\le j\le m(k,N,e)}c_{j,N,e}^kS_{k-j_1-1}^{(1+j_0)}(|D|/\delta,N,P_{j_1,k-j_1})\;,$$
where $N_2=N/2$ if $N$ is even and $N_2=N$ if $N$ is odd, $\delta=4$ if
$4\mid D$ and $\delta=1$ otherwise, we write $j=2j_1+j_0$ with $j_0\in\{0,1\}$,
upper bounds for $m(k,N,e)$ will be given below, and where we must assume
$e\ne-1$ if $N=6$ and on the contrary $e=-1$ if $N=7$.

Upper bounds for $m(k,N,e)$ are given for $e=-1$, $0$, and $1$ as follows:
$((k-1)/4,(k-1)/3,(k-3)/4)$ for $N=1$, $((k-1)/4,(k-1)/2,(k-3)/4)$ for $N=2$,
$((k-1)/2,(2k-1)/3,(k-1)/2)$ for $N=3$, $((3k-2)/4,k-1,(3k-5)/4)$ for $N=5$,
$(*,k-1,k-1)$ for $N=6$, and $(k-1,*,*)$ for $N=7$, where
$*$ denotes impossible cases.\end{conjecture}

For concreteness, we give the special case $k=3$, $e=1$: if $D\equiv1\pmod{8}$
is a negative fundamental discriminant, we have
\begin{align*}
L(\chi_D,-2)&=\dfrac{1}{35}S_2^{(1)}(|D|,1)=\dfrac{1}{7}S_2^{(1)}(|D|,2)\;,\\
(1-\lgs{D}{3})L(\chi_D,-2)&=-\dfrac{2}{63}(S_2^{(1)}(|D|,3)+14S_2^{(2)}(|D|,3))\;,\\
(1+\lgs{D}{5})L(\chi_D,-2)&=-\dfrac{2}{3}(S_2^{(1)}(|D|,5)+4S_2^{(2)}(|D|,5))\;,\\
(1-\lgs{D}{3})L(\chi_D,-2)&=\dfrac{1}{14}(-52S_2^{(1)}(|D|,6)+5S_0^{(1)}(|D|,6,1-3x))\;.
\end{align*}
Similarly to the case of even $k$, computing the number of terms involved
in the sums leads to the following algorithm:

\begin{enumerate}
\item When $D\equiv0\pmod4$: if $3\mid D$ use $N=6$, otherwise if $5\mid D$
  use $N=5$, otherwise if $D\equiv2\pmod3$ use $N=6$, otherwise if
  $D\equiv\pm1\pmod5$ use $N=5$, otherwise use $N=2$.
\item When $D\equiv1\pmod4$: if $7\mid D$ and $D\equiv5\pmod8$ use $N=7$,
  otherwise if $3\mid D$ and $D\equiv 1\pmod8$ use $N=6$,
  otherwise if $5\mid D$ use $N=5$, otherwise if $D\equiv5\pmod 8$ and
  $D\equiv 3,4,6\pmod7$ use $N=7$, otherwise if $D\equiv2\pmod3$ and
  $D\equiv1\pmod8$ use $N=6$, otherwise if $3\mid D$ (hence $D\equiv5\pmod8$)
  use $N=3$, otherwise if $D\equiv\pm1\pmod5$ use $N=5$, otherwise use $N=2$.
\end{enumerate}

As in the case of $k$ even, care must be taken in choosing $N>1$ since
the size of the linear system to be solved in order to compute the universal
coefficients $c_{j,N,e}^k$ is larger, so the above algorithm is valid
only if this time is negligible.

\bigskip

We thus give a table of timings using this algorithm. Note that $-10^j-3$
is usually (but not always) slower than $-10^j-4$ since in the latter case
the sums involve $|D|/4$ instead of $|D|$, and that a lot depends
on divisibilities by $3$, $5$, and $7$, so the tables are only an indication:

\bigskip

\centerline{
\begin{tabular}{|r||r|r|r|r|r|r|r|r|r|r|}
\hline
$D\dd k$ & $1$ & $3$   & $5$  & $7$  & $9$  & $11$ & $13$ & $15$ & $17$ & $19$ \\
\hline\hline
$-10^{10}-4$ & $0.05$ & $0.06$ & $0.06$ & $0.07$ & $0.08$ & $0.09$ & $0.10$ & $0.10$ & $0.12$ & $0.14$ \\
$-10^{10}-3$ & $0.05$ & $0.05$ & $0.06$ & $0.06$ & $0.07$ & $0.07$ & $0.08$ & $0.09$ & $0.10$ & $0.11$ \\
$-10^{11}-4$ & $0.25$ & $0.27$ & $0.28$ & $0.31$ & $0.33$ & $0.36$ & $0.40$ & $0.44$ & $0.47$ & $0.52$ \\
$-10^{11}-3$ & $0.50$ & $0.53$ & $0.56$ & $0.60$ & $0.64$ & $0.68$ & $0.73$ & $0.79$ & $0.85$ & $0.93$ \\
$-10^{12}-4$ & $1.36$ & $1.41$ & $1.47$ & $1.55$ & $1.64$ & $1.74$ & $1.86$ & $1.99$ & $2.12$ & $2.28$ \\
$-10^{12}-3$ & $2.21$ & $2.30$ & $2.40$ & $2.53$ & $2.67$ & $2.82$ & $3.00$ & $3.20$ & $3.40$ & $3.64$ \\
$-10^{13}-4$ & $6.86$ & $7.06$ & $7.28$ & $7.54$ & $7.84$ & $8.18$ & $8.55$ & $8.98$ & $9.44$ & $9.89$ \\
$-10^{13}-3$ & $35.4$ & $35.7$ & $35.9$ & $36.0$ & $36.6$ & $36.8$ & $37.0$ & $37.1$ & $37.8$ & $38.0$ \\
$-10^{14}-4$ & $34.0$ & $34.8$ & $35.5$ & $36.3$ & $37.3$ & $38.4$ & $39.7$ & $41.1$ & $42.6$ & $44.3$ \\
\hline
\end{tabular}}

\bigskip

Note that we have included the case $k=1$ which will be discussed below.

\smallskip

As we have done in the case $k$ even, for larger values of $k$ the
time to compute the coefficients dominate, so we first give a table giving
these timings:

\bigskip

\centerline{
\begin{tabular}{|r||r|r|r|r|r|r|r|r|r|r|}
\hline
$(\lgs{D}{2},N)\dd k$  & $81$   & $101$ & $151$ & $201$ & $251$ & $301$ & $351$ & $401$ & $451$ & $501$ \\
\hline\hline
  $(1,1)$ & -- & -- & $0.07$ & $0.20$ & $0.53$ & $1.22$ & $2.28$ & $3.76$ & $6.15$ & $9.37$ \\
  $(1,2)$ & -- & -- & $0.16$ & $0.48$ & $1.23$ & $3.00$ & $5.50$ & $8.83$ & $14.3$ & $21.3$ \\
  $(1,3)$ & -- & $0.10$ & $0.58$ & $1.73$ & $4.39$ & $9.31$ & $17.5$ & $30.1$ & $50.4$ & $80.6$ \\
  $(1,5)$ & $0.21$ & $0.57$ & $2.57$ & $7.43$ & $18.4$ & $37.3$ & $70.8$ & $*$ & $*$ & $*$ \\
  $(1,6)$ & -- & $0.07$ & $0.38$ & $1.44$ & $3.90$ & $10.2$ & $18.7$ & $39.4$ & $71.7$ & $*$ \\
\hline  
  $(-1,1)$ & -- & -- & $0.07$ & $0.22$ & $0.54$ & $1.27$ & $2.30$ & $3.90$ & $6.33$ & $9.70$ \\
  $(-1,2)$ & -- & -- & $0.16$ & $0.51$ & $1.23$ & $3.11$ & $5.51$ & $9.15$ & $14.4$ & $22.0$ \\
  $(-1,3)$ & -- & $0.11$ & $0.62$ & $1.82$ & $4.59$ & $9.69$ & $18.1$ & $31.4$ & $52.5$ & $64.1$ \\
  $(-1,5)$ & $0.13$ & $0.34$ & $1.53$ & $4.95$ & $10.1$ & $20.4$ & $38.5$ & $67.2$ & $110.$ & $*$ \\
  $(-1,7)$ & $0.28$ & $0.61$ & $2.73$ & $8.15$ & $20.4$ & $41.6$ & $77.4$ & $*$ & $*$ & $*$ \\
  \hline
  $(0,1)$ & -- & -- & $0.12$ & $0.39$ & $1.05$ & $1.91$ & $3.34$ & $5.80$ & $9.33$ & $14.1$ \\
  $(0,2)$ & -- & $0.07$ & $0.40$ & $1.15$ & $2.67$ & $5.76$ & $10.6$ & $19.6$ & $29.7$ & $52.8$ \\
  $(0,3)$ & $0.08$ & $0.18$ & $0.95$ & $2.82$ & $6.80$ & $14.6$ & $27.2$ & $50.6$ & $77.1$ & $*$ \\
  $(0,5)$ & $0.25$ & $0.53$ & $2.28$ & $6.78$ & $16.0$ & $33.2$ & $60.7$ & $105.$ & $*$ & $*$ \\
  $(0,6)$ & $0.29$ & $0.58$ & $2.64$ & $7.66$ & $19.0$ & $41.3$ & $75.6$ & $110.$ & $*$ & $*$ \\
\hline
\end{tabular}}

\bigskip

For future reference, we observe that the times are very roughly
\begin{align*}10^{-10}k^4(1.5,3.6,12,46.8,12.5)&\text{\quad for $e=1$,}\\
10^{-10}k^4(1.5,3.6,12,25.4,51)&\text{\quad for $e=-1$, and}\\
10^{-10}k^4(2.2,7,18,40,50)&\text{\quad for $e=0$,}\end{align*}
where as usual $e=\lgs{D}{2}$.

\smallskip

In the next table we use $N=2$ only when $|D|$ is sufficiently large, and
the corresponding timings have a $^*$; all the other timings are obtained
with $N=1$.

\bigskip

\centerline{\small
\begin{tabular}{|l||r|r|r|r|r|r|r|r|r|r|r|r|r|}
\hline
$D\dd k$ & $21$   & $41$   & $61$   & $81$   & $101$ & $151$ & $201$ & $251$ & $301$ & $351$ & $401$ & $451$ & $501$ \\
\hline\hline
$-10^{6}-20$ & -- & -- & -- & -- & -- & $0.14$ & $0.43$ & $1.11$ & $1.99$ & $3.46$ & $5.96$ & $9.57$ & $14.4$ \\
$-10^{6}-3$  & -- & -- & -- & -- & -- & $0.10$ & $0.27$ & $0.63$ & $1.41$ & $2.50$ & $4.19$ & $6.72$ & $10.2$ \\
$-10^{7}-4$  & -- & -- & -- & -- & $0.05$ & $0.19$ & $0.51$ & $1.26$ & $2.23$ & $3.82$ & $6.45$ & $10.2$ & $15.2$ \\
$-10^{7}-3$  & -- & -- & -- & -- & $0.06$ & $0.17$ & $0.41$ & $0.87$ & $1.78$ & $3.04$ & $4.95$ & $7.69$ & $11.5$ \\
$-10^{8}-20$ & -- & -- & $0.05$ & $0.08$ & $0.13$ & $0.36$ & $0.85$ & $1.85$ & $3.17$ & $5.18$ & $8.32$ & $12.7$ & $18.5$ \\
$-10^{8}-3$  & -- & $0.05$ & $0.08$ & $0.13$ & $0.18$ & $0.44$ & $0.99$ & $1.84$ & $3.29$ & $5.24$ & $8.07$ & $11.9$ & $16.9$ \\
$-10^{9}-20$ & $0.03^*$ & $0.07^*$ & $0.13^*$ & $0.23^*$ & $0.37^*$ & $0.98$ & $2.07$ & $3.91$ & $6.49$ & $10.0$ & $15.0$ & $21.8$ & $30.1$ \\
$-10^{9}-19$ & $0.08^*$ & $0.15^*$ & $0.26^*$ & $0.43^*$ & $0.62$ & $1.44$ & $3.00$ & $5.26$ & $8.60$ & $13.0$ & $19.0$ & $26.3$ & $35.5$ \\
$-10^{10}-4$ & $0.12^*$ & $0.23^*$ & $0.42^*$ & $0.70^*$ & $1.10^*$ & $2.88^*$ & $6.13^*$ & $11.2^*$ & $18.5$ & $27.5$ & $39.2$ & $54.4$ & $72.2$ \\
$-10^{10}-3$ & $0.37^*$ & $0.60^*$ & $0.97^*$ & $1.52^*$ & $2.34^*$ & $5.35$ & $10.7$ & $18.2$ & $28.6$ & $41.7$ & $59.5$ & $80.2$ & $105.$ \\
$-10^{11}-4$ & $0.50^*$ & $0.88^*$ & $1.52^*$ & $2.45^*$ & $3.80^*$ & $9.40^*$ & $19.4^*$ & $34.1^*$ & $55.6^*$ & $82.7^*$ & $*$ & $*$ & $*$ \\
$-10^{11}-3$ & $1.62^*$ & $2.41^*$ & $3.72^*$ & $5.66^*$ & $8.54^*$ & $20.6^*$ & $42.0^*$ & $72.5^*$ & $*$ & $*$ & $*$ & $*$ & $*$ \\
$-10^{12}-4$ & $2.30^*$ & $3.61^*$ & $5.78^*$ & $9.00^*$ & $13.7^*$ & $33.1^*$ & $67.2^*$ & $117.^*$ & $*$ & $*$ & $*$ & $*$ & $*$ \\
$-10^{12}-3$ & $6.65^*$ & $9.46^*$ & $14.1^*$ & $21.1^*$ & $31.6^*$ & $75.5^*$ & $*$ & $*$ & $*$ & $*$ & $*$ & $*$ & $*$ \\
$-10^{13}-4$ & $10.5^*$ & $15.1^*$ & $22.7^*$ & $34.1^*$ & $50.7^*$ & $*$ & $*$ & $*$ & $*$ & $*$ & $*$ & $*$ & $*$ \\
$-10^{13}-3$ & $40.3^*$ & $49.3^*$ & $64.8^*$ & $88.3^*$ & $*$ & $*$ & $*$ & $*$ & $*$ & $*$ & $*$ & $*$ & $*$ \\
$-10^{14}-4$ & $48.7^*$ & $64.2^*$ & $90.4^*$ & $*$ & $*$ & $*$ & $*$ & $*$ & $*$ & $*$ & $*$ & $*$ & $*$ \\
\hline
\end{tabular}}

\bigskip

\section{The Case $k=1$}

In this brief section, we consider the case $k=1$, i.e., the problem
of computing $L(\chi,0)$ for an odd character $\chi$. Of course the
Bernoulli method as well as the approximate functional equation are
still applicable in general. But in the case $\chi=\chi_D$ with
$D<0$ there are still methods coming from modular forms. Note that
in that case for $D<-4$ we have $L(\chi_D,0)=h(D)$ which can therefore be
computed using subexponential algorithms, but it is still interesting
to look at modular-type formulas. Note that $\H_1$ is not quite but almost
a modular form of weight $3/2$, so it is not surprising that the method
given above also works for $k=1$.

For instance, we have the following result, where we refer to
Definition \ref{defsigchi} for the definition of $S_0^{(1)}$
(note that $S_0^{(2)}=S_0^{(1)}$):

\begin{proposition} Let $D$ be a negative fundamental discriminant $D$.
  \begin{enumerate}\item Set $e=\lgs{D}{2}$. We have
  $$\dfrac{S_0^{(1)}(|D|,N)}{L(\chi_D,0)}=\begin{cases}
    3(1-e)&\text{\quad when $N=1$ and $N=2$\;,}\\
    (1-\lgs{D}{3})(5-e)/2&\text{\quad when $N=3$\;,}\\
    (1+\lgs{D}{5})(1-e)/2&\text{\quad when $N=5$\;,}\\
    (1-\lgs{D}{3})(1+e)/2&\text{\quad when $N=6$\;,}\\
    (1-\lgs{D}{7})&\text{\quad when $N=7$ and $e=-1$\;.}
    \end{cases}$$
\item
  If $4\mid D$, we also have
  $$\dfrac{S_0^{(1)}(|D|/4,N)}{L(\chi_D,0)}=\begin{cases}
    3&\text{\quad when $N=1$\;,}\\
    1&\text{\quad when $N=2$\;,}\\
    (1-\lgs{D}{3})/2&\text{\quad when $N=3$ and $N=6$\;,}\\
    (1+\lgs{D}{5})/2&\text{\quad when $N=5$\;.}
  \end{cases}$$
\end{enumerate}\end{proposition}

In particular, Conjecture \ref{conjodd} is valid for $k=1$ with $m(1,N,e)=0$,
$c_{0,N,e}^1=2/(3(1-e))$, $2/(3(1-e))$, $2/(5-e)$, $2/(1-e)$, $2/(1+e)$, and
$1$ when $\delta=1$ for $N=1$, $2$, $3$, $5$, $6$, and $7$ respectively,
and $c_{0,N,0}^1=2/3$, $2$, $2$, $2$, and $2$ when $\delta=4$ and
$N=1$, $2$, $3$, $5$, and $6$ respectively.

Since we can efficiently compute $L(\chi_D,0)$ by using class numbers this
result has no computational advantage, but is simply given to show that
the formulas that we obtained above for $k\ge3$ odd have analogs for $k=1$.

\end{document}